# Determination of some generalised Euler sums involving the digamma function

Donal F. Connon

dconnon@btopenworld.com

11 March 2008

**Abstract**

This paper evaluates some generalised Euler sums involving the digamma function.

**Introduction**

In response to a letter from Goldbach in 1742 (see [2, p.253]), Euler commenced his investigation of infinite series of the form

$$\sum_{n=1}^{\infty} \frac{H_n^{(p)}}{n^q}$$

where $H_n^{(p)}$ is the generalised harmonic number defined by $H_n^{(p)} = \sum_{k=1}^{n} \frac{1}{k^p}$.

In this short paper we consider related series of the form

$$\sum_{n=0}^{\infty} \frac{\psi\left(n+\frac{1}{2}\right)}{(2n+1)^2} = -\frac{1}{8}[\gamma\pi^2 + 7\varsigma(3)]$$

$$\sum_{n=0}^{\infty} \frac{\psi\left(n+\frac{1}{2}\right)}{(2n+1)^4} = -\frac{1}{96}\left[3\pi^2\varsigma(3) + \pi^4\gamma + 93\varsigma(5)\right]$$

where $\psi(x)$ is the digamma function which is the logarithmic derivative of the gamma function $\psi(x) = \frac{d}{dx}\log\Gamma(x)$ and $\varsigma(s)$ is the Riemann zeta function.

We note that [12, p.20]

$$\psi\left(n+\frac{1}{2}\right) = -\gamma - 2\log 2 + 2\sum_{k=1}^{n}\frac{1}{2k-1} = \psi\left(\frac{1}{2}\right) + 2\sum_{k=1}^{n}\frac{1}{2k-1}$$

where $\gamma$ is Euler's constant. References to some of the extensive work undertaken on Euler sums may be found in [4].

Using contour integration, de Doelder [4b] has considered series of the form

$$\sum_{n=1}^{\infty}\frac{(-1)^{n-1}}{n^2}\left[\psi\left(n+\frac{1}{2}\right)-\psi\left(\frac{1}{2}\right)\right]=2\pi G-\frac{7}{2}\varsigma(3)$$

as corrected in equation (6.30r) of [4a]. In the above equation $G$ is Catalan's constant defined by

$$G=\sum_{n=0}^{\infty}\frac{(-1)^n}{(2n+1)^2}$$

**Proof**

In 1985 Berndt [1] gave an elementary proof of Kummer's Fourier series expansion for $\log\Gamma(x)$

(1) $$\log\Gamma(x)=\frac{1}{2}\log\pi-\frac{1}{2}\log\sin\pi x+\sum_{k=1}^{\infty}\frac{(\gamma+\log 2\pi k)\sin 2k\pi x}{\pi k}\ ,\ 0<x<1$$

(this formula was originally derived by Kummer in 1847 [9]). Reference to the well-known Fourier series [13, p.148]

$$\log[2\sin\pi x]=-\sum_{k=1}^{\infty}\frac{\cos 2k\pi x}{k}\ ,\ 0<x<1$$

confirms that (1) is properly described as a Fourier series expansion for $\log\Gamma(x)$

(1a) $$\log\Gamma(x)=\frac{1}{2}\log(2\pi)+\frac{1}{2}\sum_{k=1}^{\infty}\frac{\cos 2k\pi x}{k}+\sum_{k=1}^{\infty}\frac{(\gamma+\log 2\pi k)\sin 2k\pi x}{\pi k}$$

We now multiply (1a) across by $\sin(2n+1)\pi x$ and integrate term by term. We have the familiar trigonometric formula

$$2\cos 2k\pi x\sin(2n+1)\pi x=\sin(2k+2n+1)\pi x-\sin(2k-2n-1)\pi x$$

and integration results in

$$2\int_0^1\cos 2k\pi x\sin(2n+1)\pi x\,dx=-\frac{\cos(2k+2n+1)\pi x}{(2k+2n+1)\pi}+\frac{\cos(2k-2n-1)\pi x}{(2k-2n-1)\pi}\bigg|_0^1$$

$$=\frac{2}{(2k+2n+1)\pi}-\frac{2}{(2k-2n-1)\pi}$$

We then have

$$S_1 = \frac{1}{2}\sum_{k=1}^{\infty}\frac{1}{k}\int_0^1 \cos 2k\pi x \sin(2n+1)\pi x\, dx = \frac{1}{2\pi}\sum_{k=1}^{\infty}\left[\frac{1}{k(2k+2n+1)} - \frac{1}{k(2k-2n-1)}\right]$$

By partial fractions and letting $a = 2n+1$ we see that

$$\frac{1}{(2k+a)k} = \frac{1}{a}\left(\frac{1}{k} - \frac{2}{2k+a}\right)$$

and therefore

$$\frac{1}{(2k+a)k} - \frac{1}{(2k-a)k} = \frac{1}{a}\left(\frac{1}{k} - \frac{2}{2k+a}\right) + \frac{1}{a}\left(\frac{1}{k} - \frac{2}{2k-a}\right)$$

$$= \frac{1}{a}\left[\frac{1}{k} - \frac{1}{k+(a/2)} + \frac{1}{k} - \frac{1}{k-(a/2)}\right]$$

Then, using the well known expression for the digamma function [14, p.108]

$$\psi(x) = -\gamma - \frac{1}{x} + \sum_{k=1}^{\infty}\left(\frac{1}{k} - \frac{1}{k+x}\right)$$

we see that

$$S_1 = \frac{1}{2\pi a}\left[\psi(a/2) + \psi(-a/2) + 2\gamma\right]$$

We have the well known formula [12, p.14]

$$\psi(-x) = \psi(x) + \pi \cot \pi x + \frac{1}{x}$$

and hence we get

$$S_1 = \frac{1}{2\pi a}\left[2\psi(a/2) + \pi \cot(\pi a/2) + \frac{2}{a} + 2\gamma\right]$$

Since $a = 2n+1$ this becomes

$$S_1 = \frac{1}{(2n+1)\pi}\left[\psi\left(n+\frac{1}{2}\right) + \frac{1}{2n+1} + \gamma\right]$$

and using the identity [12, p.20]

$$\psi\left(n+\frac{1}{2}\right) = -\gamma - 2\log 2 + 2\sum_{k=1}^{n}\frac{1}{2k-1}$$

we obtain

$$S_1 = \frac{1}{(2n+1)\pi}\left[-2\log 2 + 2\sum_{k=1}^{n}\frac{1}{2k-1} + \frac{1}{2n+1}\right]$$

It is readily seen that

$$\frac{1}{2}\int_0^1 \log(2\pi)\sin(2n+1)\pi x\,dx = \frac{\log(2\pi)}{(2n+1)\pi}$$

and since $\int_0^1 \sin(2n+1)\pi x \cos 2k\pi x = 0$ we deduce that

(2) $$\int_0^1 \log\Gamma(x)\sin(2n+1)\pi x\,dx = \frac{1}{(2n+1)\pi}\left[\log\left(\frac{\pi}{2}\right) + \frac{1}{2n+1} + 2\sum_{k=1}^{n}\frac{1}{2k-1}\right]$$

As noted by Kölbig [8], this integral was incorrectly reported in Nielsen's book [11, p.203]. The result was corrected in the fifth edition of Gradshteyn and Ryzhik [6, No.6.443 2]).

The integral in (2) may also be written as

(2a) $$\int_0^1 \log\Gamma(x)\sin(2n+1)\pi x\,dx = \frac{1}{(2n+1)\pi}\left[\log(2\pi) + \gamma + \frac{1}{2n+1} + \psi\left(n+\frac{1}{2}\right)\right]$$

We then multiply (2a) by $\dfrac{1}{2n+1}$ and make the summation

$$S_2 = \sum_{n=0}^{\infty}\frac{1}{2n+1}\int_0^1 \log\Gamma(x)\sin(2n+1)\pi x\,dx$$

$$= \frac{\log(2\pi)+\gamma}{\pi}\sum_{n=0}^{\infty}\frac{1}{(2n+1)^2}+\frac{1}{\pi}\sum_{n=0}^{\infty}\frac{1}{(2n+1)^3}+\frac{1}{\pi}\sum_{n=0}^{\infty}\frac{\psi\left(n+\frac{1}{2}\right)}{(2n+1)^2}$$

which may also be written as

$$S_2 = \sum_{n=0}^{\infty}\frac{1}{2n+1}\int_0^1 \log\Gamma(x)\sin(2n+1)\pi x\,dx$$

$$= \frac{1}{\pi}\log\pi\sum_{n=0}^{\infty}\frac{1}{(2n+1)^2}+\frac{1}{\pi}\sum_{n=0}^{\infty}\frac{1}{(2n+1)^3}+\frac{2}{\pi}\sum_{n=0}^{\infty}\frac{1}{(2n+1)^2}\sum_{k=1}^{n}\frac{1}{2k-1}$$

It is well known that

$$\sum_{n=0}^{\infty}\frac{1}{(2n+1)^s} = (1-2^{-s})\varsigma(s)$$

and we therefore have

$$\sum_{n=0}^{\infty}\frac{1}{(2n+1)^2} = \frac{\pi^2}{8} \qquad \sum_{n=0}^{\infty}\frac{1}{(2n+1)^3} = \frac{7}{8}\varsigma(3)$$

Hence we obtain

(3) $$S_2 = \frac{\pi}{8}[\log(2\pi)+\gamma]+\frac{7\varsigma(3)}{8\pi}+\frac{1}{\pi}\sum_{n=0}^{\infty}\frac{\psi\left(n+\frac{1}{2}\right)}{(2n+1)^2}$$

Assuming we may validly change the order of summation and integration, we have

$$\sum_{n=0}^{\infty}\frac{1}{2n+1}\int_0^1 \log\Gamma(x)\sin(2n+1)\pi x\,dx = \int_0^1 \log\Gamma(x)\sum_{n=0}^{\infty}\frac{\sin(2n+1)\pi x}{2n+1}\,dx$$

For $0 < x < 1$ we have the Fourier expansion [13, p.149]

$$\sum_{n=0}^{\infty}\frac{\sin(2n+1)\pi x}{2n+1} = \frac{\pi}{4}$$

and therefore we obtain for the left-hand side

(4) $$S_2 = \frac{\pi}{4}\int_0^1 \log\Gamma(x)\,dx = \frac{\pi}{8}\log(2\pi)$$

where we have used Raabe's integral [15, p.261] $\int_0^1 \log \Gamma(x)\,dx = \frac{1}{2}\log(2\pi)$.

Equating (3) and (4) gives us

(5) $$\sum_{n=0}^{\infty} \frac{\psi\left(n+\frac{1}{2}\right)}{(2n+1)^2} = -\frac{1}{8}[\gamma\pi^2 + 7\varsigma(3)]$$

Independent confirmation of the above identity was soon forthcoming because, most fortuitously, I also spotted it in Milgram's paper [10]; this identity was originally obtained by Jordan [7] in 1973.

Emboldened by the success of the first example, we now continue this procedure. As before, the identity in (1a) is multiplied by $\frac{1}{(2n+1)^3}$ and we make the summation

$$\sum_{n=0}^{\infty} \frac{1}{(2n+1)^3} \int_0^1 \log\Gamma(x)\sin(2n+1)\pi x\,dx =$$

$$\frac{\log(2\pi)+\gamma}{\pi}\sum_{n=0}^{\infty}\frac{1}{(2n+1)^4} + \frac{1}{\pi}\sum_{n=0}^{\infty}\frac{1}{(2n+1)^5} + \frac{1}{\pi}\sum_{n=0}^{\infty}\frac{\psi\left(n+\frac{1}{2}\right)}{(2n+1)^4}$$

(6) $$= \frac{\pi^3}{96}[\log(2\pi)+\gamma] + \frac{31\varsigma(5)}{32\pi} + \frac{1}{\pi}\sum_{n=0}^{\infty}\frac{\psi\left(n+\frac{1}{2}\right)}{(2n+1)^4}$$

Assuming we may validly change the order of summation and integration, we have

$$\sum_{n=0}^{\infty}\frac{1}{(2n+1)^3}\int_0^1 \log\Gamma(x)\sin(2n+1)\pi x\,dx = \int_0^1 \log\Gamma(x)\sum_{n=0}^{\infty}\frac{\sin(2n+1)\pi x}{(2n+1)^3}\,dx$$

We have the well-known Fourier series from Tolstov's book [13, p.149] for $0 \leq x \leq 1$

$$\sum_{n=0}^{\infty}\frac{\sin(2n+1)\pi x}{(2n+1)^3} = \frac{\pi^3}{8}(x-x^2)$$

and hence we have

$$\int_0^1 \log\Gamma(x) \sum_{n=0}^{\infty} \frac{\sin(2n+1)\pi x}{(2n+1)^3} dx = \frac{\pi^3}{8} \int_0^1 (x-x^2) \log\Gamma(x) dx$$

We note that Espinosa and Moll [5] report the following two integrals

(6a) $$\int_0^1 x \log\Gamma(x) dx = \frac{1}{4}\log(2\pi) - \log A$$

(6b) $$\int_0^1 x^2 \log\Gamma(x) dx = \frac{\varsigma(3)}{4\pi^2} + \frac{1}{6}\log(2\pi) - \log A$$

where the Glaisher-Kinkelin constant $A$ is defined by

$$\log A = \frac{1}{12} - \varsigma'(-1)$$

We then obtain

(7) $$\frac{\pi^3}{8}\int_0^1 (x-x^2)\log\Gamma(x) dx = \frac{\pi^3}{96}\log(2\pi) - \frac{\pi\varsigma(3)}{32}$$

and equating (6) and (7) we get

(8) $$\sum_{n=0}^{\infty} \frac{\psi\left(n+\frac{1}{2}\right)}{(2n+1)^4} = -\frac{1}{96}\left[3\pi^2\varsigma(3) + \pi^4\gamma + 93\varsigma(5)\right]$$

On the assumption that $\gamma$ and $\pi$ each have a weight of one and that $\varsigma(n)$ has a weight of $n$, we see that all of the terms of the above identity have a weight equal to five (and this feature of homogeneity is also present in the standard Euler sums).

In a similar way we may be able to determine an expression for $\sum_{n=0}^{\infty} \frac{\psi\left(n+\frac{1}{2}\right)}{(2n+1)^4}$ provided we can evaluate $\int_0^1 x^4 \log\Gamma(x) dx$. As a minimum, this will require the following result from [3] and [12, p.209]

$$2\int_0^z t \log\Gamma(a+t) dt = \left(\frac{1}{4} - \frac{1}{2}a + \frac{1}{2}a^2 - 2\log A\right)z + \left(\frac{1}{2}\log(2\pi) - \frac{1}{2}a + \frac{1}{4}\right)z^2$$

$$-\frac{1}{2}z^2 + z^2 \log \Gamma(z+a) - (a-1)^2 [\log \Gamma(z+a) - \log \Gamma(a)] + (2a-3)[\log G(z+a) - \log G(a)]$$

$$-2[\log \Gamma_3(z+a) - \log \Gamma_3(a)]$$

where the double and triple gamma functions (due to Barnes) are defined by

$$[\Gamma_2(1+x)]^{-1} = G(1+x) = (2\pi)^{x/2} \exp\left[-\frac{1}{2}(\gamma x^2 + x^2 + x)\right] \prod_{k=1}^{\infty} \left\{\left(1+\frac{x}{k}\right)^k \exp\left(\frac{x^2}{2k} - x\right)\right\}$$

and

$$\Gamma_3(1+x) = G_3(1+x) = \exp(c_1 x + c_2 x^2 + c_3 x^3) F(x)$$

where

$$F(x) = \prod_{k=1}^{\infty} \left\{\left(1+\frac{x}{k}\right)^{-\frac{1}{2}k(k+1)} \exp\left(\frac{1}{2}\left(1+\frac{1}{k}\right)kx - \frac{1}{4}\left(1+\frac{1}{k}\right)x^2 + \frac{1}{6k}\left(1+\frac{1}{k}\right)x^3\right)\right\}$$

and the constants are

$$c_1 = \frac{3}{8} - \frac{1}{4}\log(2\pi) - \log A, \qquad c_2 = \frac{1}{4}\left[\gamma + \log(2\pi) + \frac{1}{2}\right]$$

$$c_3 = -\frac{1}{6}\left[\gamma + \varsigma(2) + \frac{3}{2}\right]$$

In this regard, Jordan [7] has reported that

$$4^{n+1} \sum_{k=1}^{\infty} \frac{L_k}{(2k-1)^{2n}} = \varsigma(2n+1) + 4\varsigma(2n)\log 2 - 2\sum_{j=1}^{n-1} \varsigma(2j)\varsigma(2n+1-2j)$$

where $L_k = \sum_{j=1}^{k} \frac{1}{2j-1} = \frac{1}{2}\left[\psi\left(k+\frac{1}{2}\right) - \psi\left(\frac{1}{2}\right)\right]$.

Kölbig [8] also showed that

$$\int_0^1 \log\Gamma(x)\cos(2n+1)\pi x\,dx = \frac{2}{\pi^2}\left[\frac{\log(2\pi)+\gamma}{(2n+1)^2} + 2\sum_{k=1}^{\infty} \frac{\log k}{4k^2 - (2n+1)^2}\right]$$

and completing the summation we have

$$\sum_{n=1}^{\infty}\frac{1}{(2n+1)^2}\int_0^1 \log\Gamma(x)\cos(2n+1)\pi x\,dx = \frac{2[\log(2\pi)+\gamma]}{\pi^2}\sum_{n=1}^{\infty}\frac{1}{(2n+1)^4}$$

$$+\frac{4}{\pi^2}\sum_{n=1}^{\infty}\frac{1}{(2n+1)^2}\sum_{k=1}^{\infty}\frac{\log k}{4k^2-(2n+1)^2}$$

(9)
$$=\frac{[\log(2\pi)+\gamma]\pi^2}{48}$$

$$+\frac{4}{\pi^2}\sum_{n=1}^{\infty}\frac{1}{(2n+1)^2}\sum_{k=1}^{\infty}\frac{\log k}{4k^2-(2n+1)^2}$$

Assuming we may validly change the order of summation and integration, we have

$$\sum_{n=1}^{\infty}\frac{1}{(2n+1)^2}\int_0^1 \log\Gamma(x)\cos(2n+1)\pi x\,dx = \int_0^1 \log\Gamma(x)\sum_{n=1}^{\infty}\frac{\cos(2n+1)\pi x}{(2n+1)^2}\,dx$$

and we then employ the Fourier series [13, p.149] for $0 \le t \le \pi$

$$\sum_{n=1}^{\infty}\frac{\cos(2n+1)t}{(2n+1)^2} = \frac{\pi^2-2\pi t}{8}$$

to obtain

$$=\frac{\pi^2}{8}\int_0^1 \log\Gamma(x)[1-2x]\,dx$$

Using (6a) this becomes

(10)
$$=\frac{\pi^2}{8}\left[\frac{1}{4}\log(2\pi)+\log A\right]$$

Equating (9) and (10) gives us

(11)
$$\sum_{n=1}^{\infty}\frac{1}{(2n+1)^2}\sum_{k=1}^{\infty}\frac{\log k}{4k^2-(2n+1)^2} = \frac{\pi^4}{384}[\log(2\pi)+12\log A-2\gamma]$$

In passing, we also note from [11, p.203]

$$\int_0^1 \log \Gamma(x) \cos 2n\pi x \, dx = \frac{1}{4n}$$

and hence the summation gives us

$$\sum_{n=1}^{\infty} \frac{1}{n^2} \int_0^1 \log \Gamma(x) \cos 2n\pi x \, dx = \frac{1}{4}\varsigma(3)$$

Assuming we may validly change the order of summation and integration, we have

$$\sum_{n=1}^{\infty} \frac{1}{n^2} \int_0^1 \log \Gamma(x) \cos 2n\pi x \, dx = \int_0^1 \log \Gamma(x) \sum_{n=1}^{\infty} \frac{\cos 2n\pi x}{n^2} \, dx$$

and we then employ the well-known Fourier series [13, p.148]

$$\sum_{n=1}^{\infty} \frac{\cos 2n\pi x}{n^2} = \pi^2 \left( x^2 - x + \frac{1}{6} \right) = \pi^2 B_2(x)$$

and obtain the definite integral

$$\int_0^1 B_2(x) \log \Gamma(x) \, dx = \frac{\varsigma(3)}{4\pi^2}$$

which is in agreement with [5].

In a recent preprint entitled "Euler-Hurwitz series" (to be submitted to arXiv) the author has employed (5) to compute the following integral

$$\int_0^1 \frac{\sqrt{v}}{1-v} \left[ \frac{7}{8}\varsigma(3) - \frac{1}{4}\Phi\left(v, 2, \frac{1}{2}\right) \right] dv = \frac{21}{8}\varsigma(3) - \frac{3}{2}\log 2\varsigma(2)$$

where $\Phi(v, s, a)$ is the Hurwitz-Lerch zeta function [12, p.121] defined by

$$\Phi(v, s, a) = \sum_{n=0}^{\infty} \frac{v^n}{(n+a)^s}$$

Series involving the digamma function were extensively explored by Coffey [3a] in 2005. His paper included a number of series of the form

$$\sum_{n=1}^{\infty}\frac{1}{n^2}\sum_{r=0}^{q-1}\left[\gamma+\psi\left(n+\frac{r}{q}\right)\right]=q\left[\left(q^2+\frac{1}{2q}\right)\varsigma(3)+\pi\sum_{j=1}^{q-1}\text{Cl}_2\left(\frac{2\pi j}{q}\right)\right]-q\varsigma(2)\log q$$

where the Clausen integral [12, p.111] is defined by

$$\text{Cl}_2(t)=-\int_0^t\log\left[2\sin\frac{x}{2}\right]dx=\sum_{n=1}^{\infty}\frac{\sin nt}{n^2}\quad,0\le x\le 2\pi$$

However, it is not immediately clear to the author how equations (5) and (8) above may be derived directly from Coffey's paper [3a].

**REFERENCES**

[1] B.C. Berndt, The Gamma Function and the Hurwitz Zeta Function.
    Amer. Math. Monthly, 92,126-130, 1985.

[2] B.C. Berndt, Ramanujan's Notebooks. Part I, Springer-Verlag, 1989.

[3] J. Choi and H.M. Srivastava, Certain classes of series associated with the Zeta
    function and multiple gamma functions.
    J. Comput. Appl. Math., 118, 87-109, 2000.

[3a] M.W. Coffey, On one-dimensional digamma and polygamma series related to
     the evaluation of Feynman diagrams.
     J. Comput. Appl. Math, 183, 84-100, 2005. math-ph/0505051 [abs, ps, pdf, other]

[4] D.F. Connon, Some series and integrals involving the Riemann zeta function,
    binomial coefficients and the harmonic numbers. Volume I, 2007.
    arXiv:0710.4022 [pdf]

[4a] D.F. Connon, Some series and integrals involving the Riemann zeta function,
     binomial coefficients and the harmonic numbers. Volume V, 2007.
     arXiv:0710.4047 [pdf]

[4b] P.J. de Doelder, On some series containing $\psi(x)-\psi(y)$ and $(\psi(x)-\psi(y))^2$ for
     certain values of $x$ and $y$. J. Comput. Appl. Math. 37, 125-141, 1991.

[5] O. Espinosa and V.H. Moll, On some integrals involving the Hurwitz zeta
    function: Part I. The Ramanujan Journal, 6,150-188, 2002.
    http://arxiv.org/abs/math.CA/0012078

[6] I.S. Gradshteyn and I.M. Ryzhik, Tables of Integrals, Series and Products.
    Sixth Ed., Academic Press, 2000.

Errata for Sixth Edition http://www.mathtable.com/errata/gr6_errata.pdf

[7] P.F. Jordan, Infinite sums of Psi functions.
Bull. Amer. Math. Soc. 79 (1973), 681-683.
http://www.ams.org/bull/1973-79-04/S0002-9904-1973-13258-6/

[8] K.S. Kölbig, On three trigonometric integrals of $\log \Gamma(x)$ or its derivatives.
CERN/Computing and Networks Division; CN/94/7, May 1994.
J. Comput. Appl. Math. 54 (1994) 129-131.
http://doc.cern.ch/tmp/convert_P00023218.pdf

[9] E.E. Kummer, Beitrag zur Theorie der Function $\Gamma(x) = \int_0^\infty e^{-v} v^{x-1} dv$.
J. Reine Angew. Math., 35, 1-4, 1847.
http://www.digizeitschriften.de/index.php?id=132&L=2

[10] M. Milgram, On Some Sums of Digamma and Polygamma functions.
arXiv:math/0406338 [pdf] ,2004.

[11] N. Nielsen, Die Gammafunktion. Chelsea Publishing Company, Bronx and New York, 1965.

[12] H.M. Srivastava and J. Choi, Series Associated with the Zeta and Related Functions. Kluwer Academic Publishers, Dordrecht, the Netherlands, 2001.

[13] G.P. Tolstov, Fourier Series. (Translated from the Russian by R.A. Silverman) Dover Publications Inc, New York, 1976.

[14] Z.X. Wang and D.R. Guo, Special Functions.
World Scientific Publishing Co. Pte. Ltd, Singapore, 1989.

[15] E.T. Whittaker and G.N. Watson, A Course of Modern Analysis: An Introduction to the General Theory of Infinite Processes and of Analytic Functions; With an Account of the Principal Transcendental Functions. Fourth Ed., Cambridge University Press, Cambridge, London and New York, 1963.